\documentclass{article}

\usepackage{arxiv}

\usepackage[utf8]{inputenc} 
\usepackage[T1]{fontenc}    
\usepackage{hyperref}       
\usepackage{url}            
\usepackage{booktabs}       
\usepackage{amsfonts}       
\usepackage{nicefrac}       
\usepackage{microtype}      
\usepackage{graphicx}

\usepackage{doi}

\usepackage{graphicx}
\usepackage{mathtools}
\usepackage{amsfonts}
\usepackage[style=numeric, giveninits=true]{biblatex}
\bibliography{AC_KK_Heegard.bib}

\newtheorem{theorem}{Theorem}[section]

\title{Each Thurston geometry admits a Heegaard genus two example.}

\author{ \href{https://orcid.org/
0000-0003-0877-0731}{\includegraphics[scale=0.06]{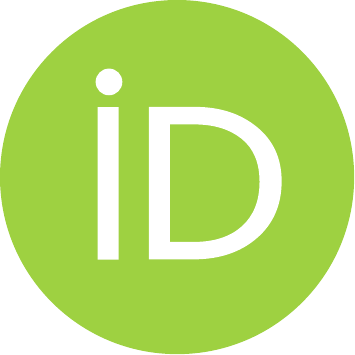}\hspace{1mm}Andrzej Czarnecki} \\
	Jagiellonian University\\
	Faculty of Mathematics and Computer Science\\
	\L{}ojasiewicza 6 \\
	30-348, Krak\'{o}w, Poland\\
	\texttt{andrzejczarnecki01@gmail.com} \\
	\And
	\href{https://orcid.org/0000-0002-9059-592X}{\includegraphics[scale=0.06]{orcid.pdf}\hspace{1mm}Katarzyna Krawiec} \\
    University of Warsaw\\
    Faculty of Mathematics, Informatics and Mechanics\\
    Banacha 2\\
    02-097, Warszawa, Poland\\
   \texttt{kkrawiec97@gmail.com}
}

\renewcommand{\shorttitle}{Each geometry admits a genus two example.}

\hypersetup{
pdftitle={A template for the arxiv style},
pdfsubject={math.AT, math.GT},
pdfauthor={Andrzej Czarnecki, Katarzyna Krawiec},
pdfkeywords={Geometrization theorem, Heegaard diagrams},
}

\begin{document}
\maketitle

\begin{abstract}
	We collect examples of genus two Heegaard diagrams for compact 3-manifolds admitting each of the eight Thurston geometries.
\end{abstract}

\keywords{Geometrization theorem \and Heegaard diagrams}

\section*{Introduction}

This paper collects the explicit examples of genus two diagrams of closed 3-manifolds of each of the Thurston geometries. While the first example, the Poincar\'e sphere, is well known and included only for completeness, the subsequent cases are less trivial. Even the simplest diagram in the paper, of $\mathbb{RP}^3\sharp\mathbb{RP}^3$, relies on the fact that that manifold is the unique non-prime manifold that admits geometrization. Most importantly however, some are of note since their genus was known, but computed non-constructively and we could find no diagrams in literature. We also could not find references to the observation made in the title.

Note that of course any 3-manifold of genus one is very restricted in terms of geometry (being either the sphere, a lens space, or $\mathbb{S}^2\times\mathbb{S}^1$). Thus the examples below show that as soon as the genus raises to two, these restrictions are lifted completely.

To ensure that our diagrams indeed present the desired manifolds, we will use the following consequence of Geometrization Theorem, and several other classical results (for the complete discussion, cf. \cite{friedl} and references therein).
\begin{theorem}[Theorem 2.1.2 in \cite{friedl}]\label{maszynka}
Let $N$ and $N'$ be closed, orientable, prime 3-manifolds and let $\phi:\pi_1(N)\rightarrow\pi_1(N')$ be an isomorphism.
\begin{enumerate}
    \item If $N$ and $N'$ are not lens spaces, then $N$ and $N'$ are homeomorphic.
    \item If $N$ and $N'$ are not spherical, then there exists a homeomorphism which induces $\phi$.
\end{enumerate}
\end{theorem}
We will treat the Heegaard diagrams as in \cite{singer}, four discs on the boundary of a closed 3-ball, $a$, $a^{-1}$, $b$, and $b^{-1}$ being pairwise identified (with appropriate orientation flip), constituting one handlebody of a Heegaard splitting, and two curves being the images of respective circles in the second handlebody (and prescribing a homeomorphism of their boundaries). This will be equivalent to prescribing two generators and two relations in the fundamental group of the manifold decomposed: if we are able to draw a diagram encoding the group's presentation, the theorem above ensures that it is indeed a diagram of the initial manifold. At this point it is important to note that not every presentation of a fundamental group is realizable in the diagrammatic fashion described above. In particular, the rank of the fundamental group is only bounded from above, not a priori equal to the Heegaard genus of the manifold. The question was first posed by Waldhausen \cite{Waldhausen} and currently it is well-known that there are counterexamples that are both Seifert-fibered \cite{Boileau-Zieschang} and hyperbolic \cite{jak}.

In what follows, we found easier to draw diagrams for presentations of groups on two generators $a$, $b$ acting on themselves on the right. This is of course only a cosmetic choice, but the reader should keep in mind that the relations are to be read from left to right. Each relation starts from the number in a bold circle, and the relation given by the dashed line is always given second. The familiar example immediately below serves as an introduction.

\section{Spherical geometry: Poincar\'e sphere}

\begin{figure}[h!]
    \centering
    \includegraphics[width=0.7\textwidth]{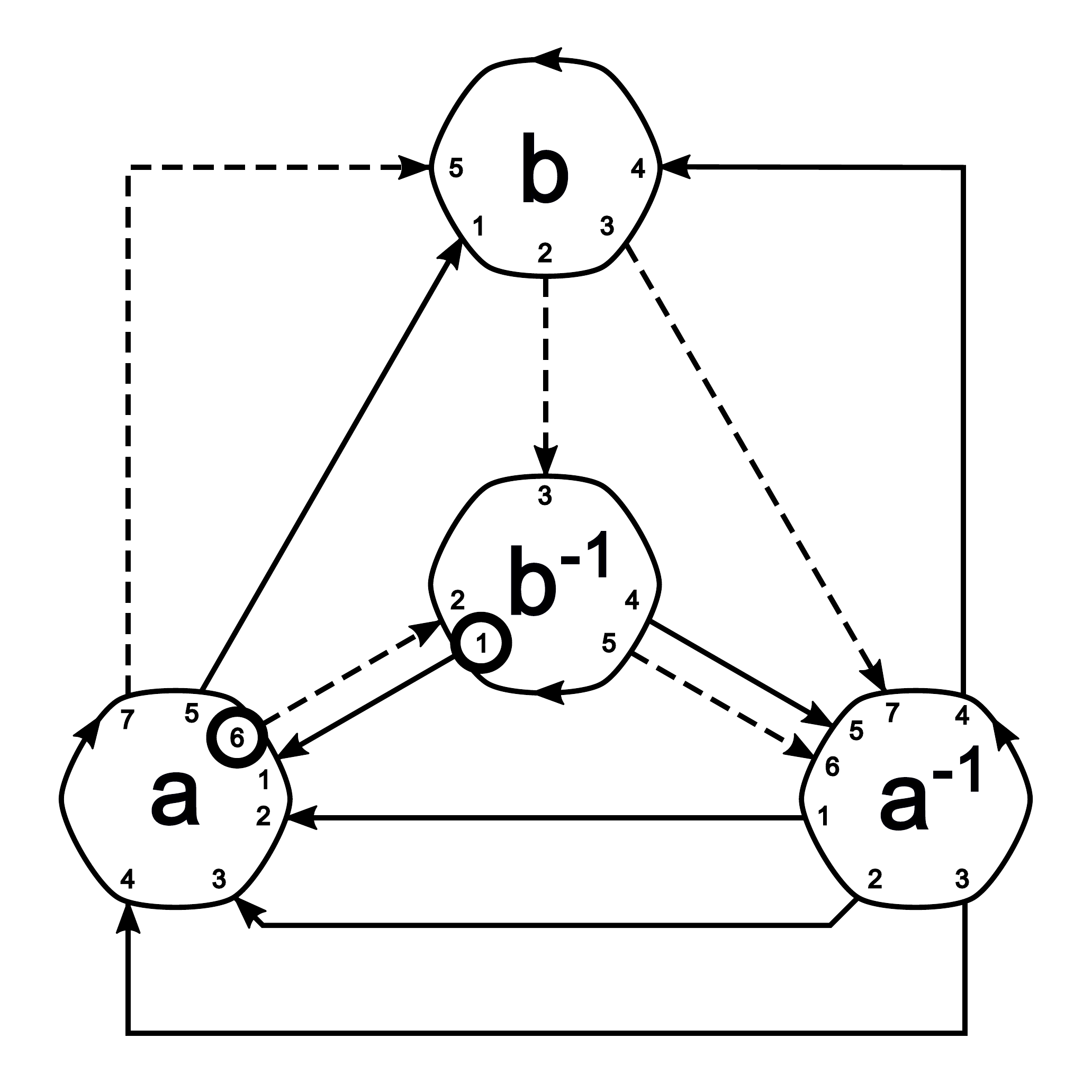}
    \caption{The two relations are $a^4ba^{-1}b=1$ and $b^{-2}a^{-1}ba^{-1}=1$.}
    \label{fig:poincare}
\end{figure}

This is the ``straightened'', but otherwise well-known diagram for the Poincar\'e sphere as originally defined in \cite{dlaczego} and afterwards repeated in many places. The presentation of the fundamental group is easily seen to be of the binary icosahedral group, which easily follows from  its usual presentation $2I=\left\langle s, t \,\middle|\,(st)^2=s^3, (st)^2=t^5\right\rangle$.

\newpage

\section{\texorpdfstring{$\mathbb{S}^2\times\mathbb{R}$}{S2xR} geometry: \texorpdfstring{$\mathbb{RP}^3\sharp\mathbb{RP}^3$}{RP3 RP3}}

\begin{figure}[h!]
    \centering
    \includegraphics[width=0.7\textwidth]{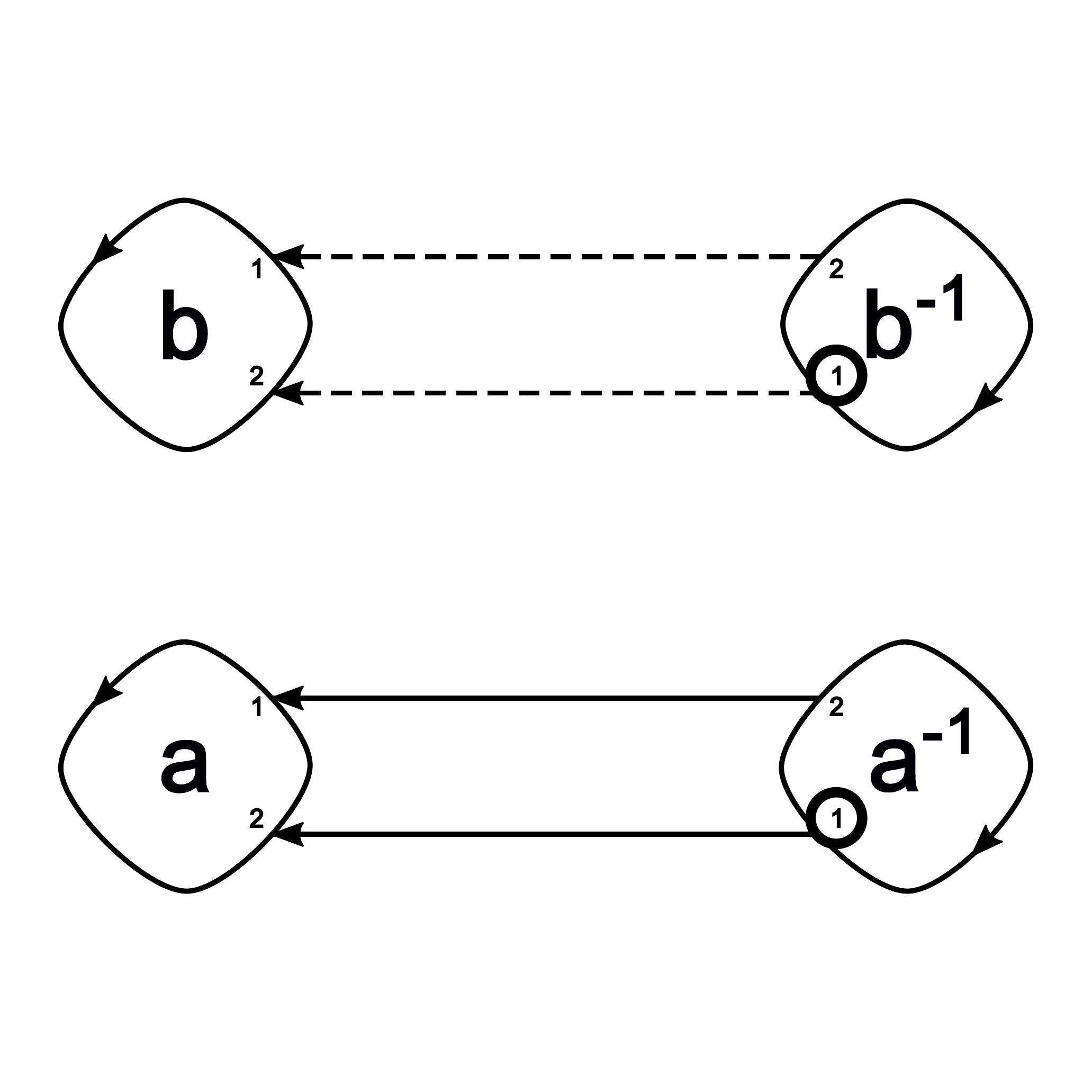}
    \caption{The two relations are $a^2=1$ and $b^2=1$.}
    \label{fig:2rp3}
\end{figure}

It is well known that there are only two orientable closed 3-manifolds with $\mathbb{S}^2\times\mathbb{R}$ geometry, the $\mathbb{S}^2\times \mathbb{S}^1$ and $\mathbb{RP}^3\sharp\mathbb{RP}^3$. The former is of Heegard genus 1, but luckily the latter is the only non-prime 3-manifold that admits a geometric structure, and is of genus 2. The straightforward diagram is included for completeness.

\newpage

\section{Hyperbolic geometry: Weeks manifold}

\begin{figure}[h!]
    \centering
    \includegraphics[width=0.7\textwidth]{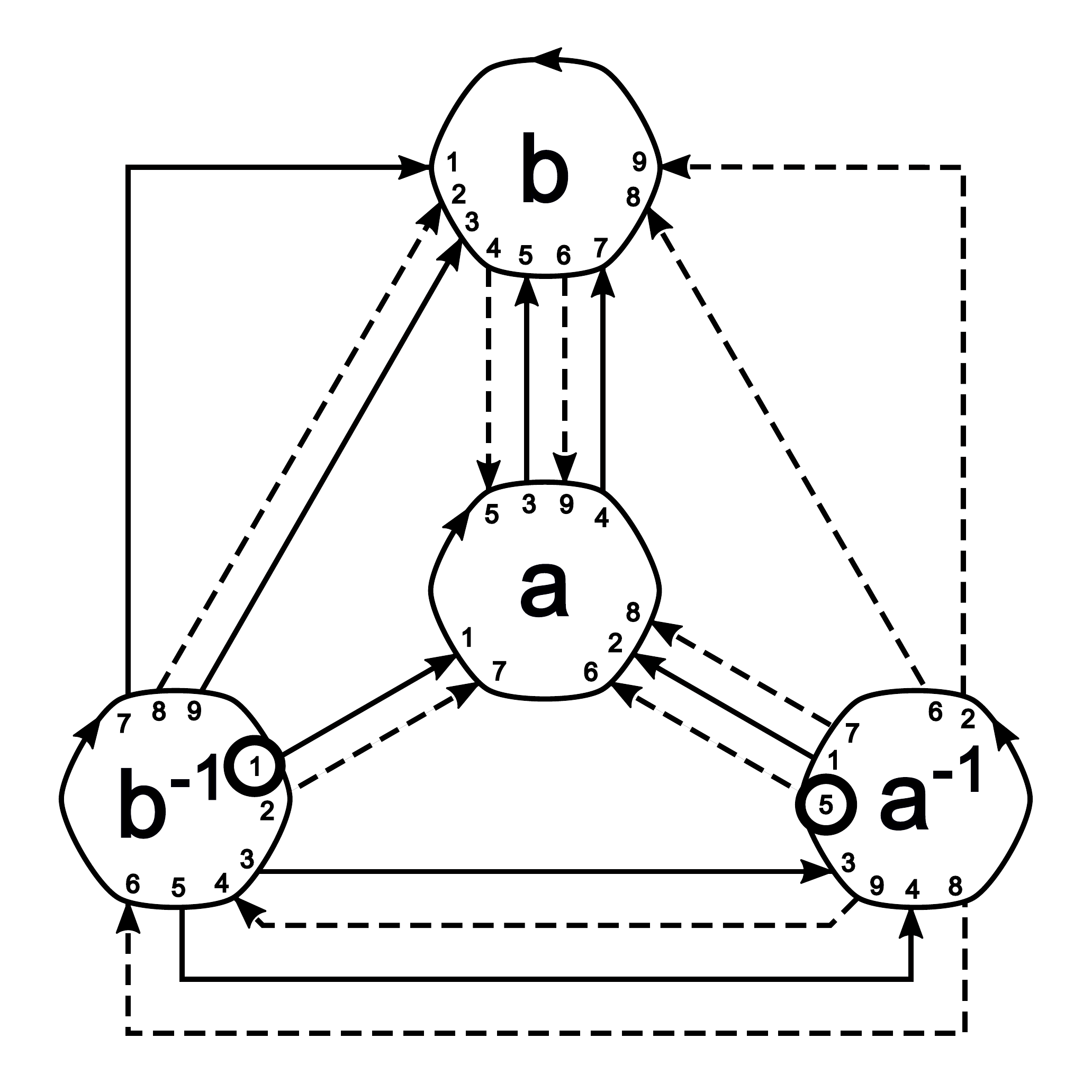}
    \caption{The two relations are $a^2b^2a^{-1}ba^{-1}b^2=1$ and $ab^2a^2b^{-1}ab^{-1}a=1$.}
    \label{fig:weeks}
\end{figure}

The Weeks manifold is notable as the smallest hyperbolic manifold (recall that the volumes of closed hyperbolic 3-manifold form a well-ordered subset of $\mathbb{R}$, and thus attain minimum, \cite{wjakimcelu}, Theorem 5.12.1). It was identified in \cite{poco} as the $(5,2)$ and $(5,1)$ Dehn surgeries on the Whitehead link, and shown to have the fundamental group presentable as above.

\newpage

\section{Flat geometry: a Fibonacci manifold}

\begin{figure}[h!]
    \centering
    \includegraphics[width=0.7\textwidth]{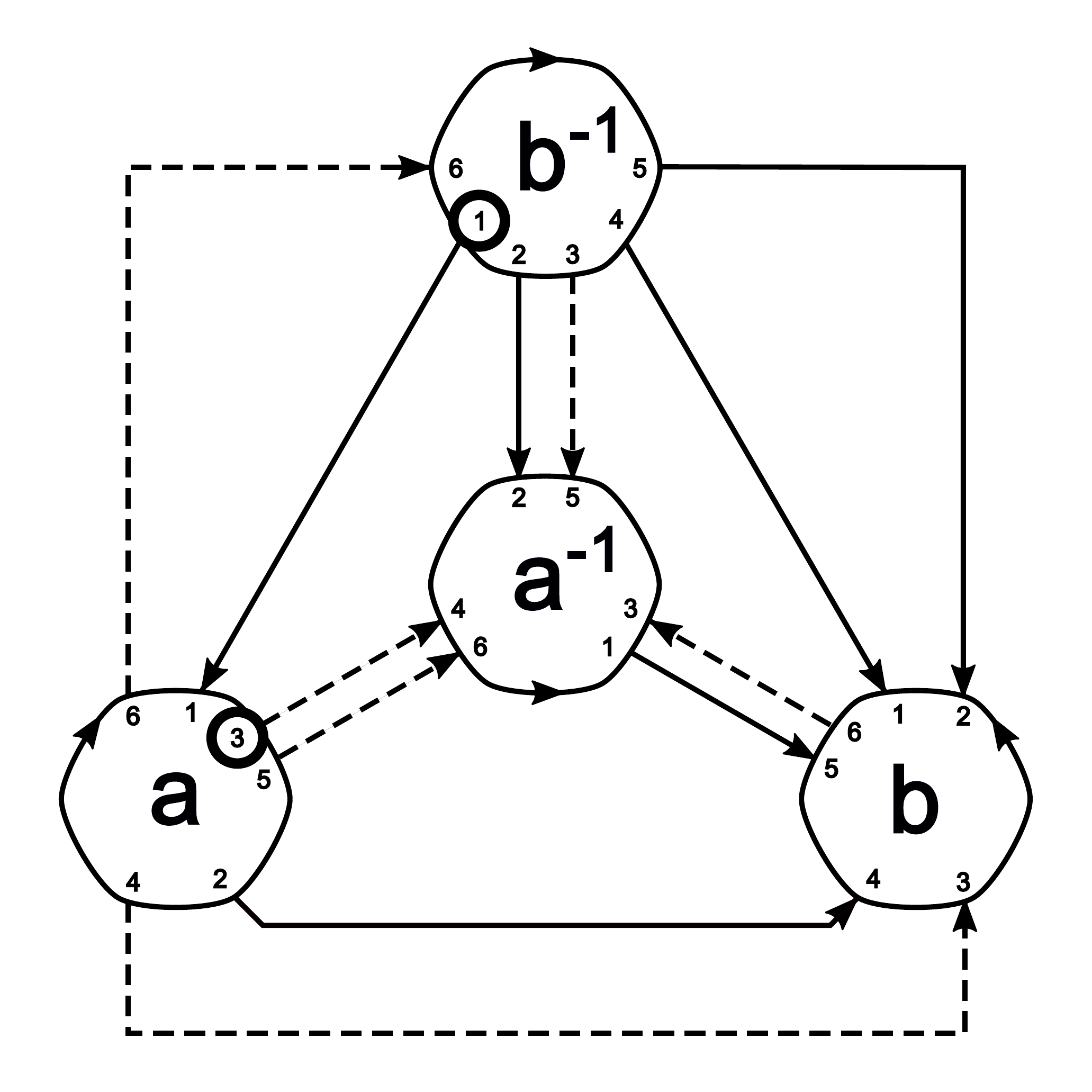}
    \caption{The two relations are $ab^2a^{-1}b^2=1$ and $a^{-1}ba^{-2}b^{-1}a^{-1}=1$.} 
    \label{fig:fibonacci}
\end{figure}

It is well known that there are only six closed and orientable flat 3-manifolds, of which only one has its first Betti number 0. This manifold, presented above, is sometimes called Hantzsche-Wendt manifold. By \cite{geometricstudyoffibonacci} (Proposition B), it is the only non-hyperbolic example of so-called Fibonacci 3-manifolds, given by their fundamental group, the Fibonacci group $F(2,2n)$ (disregarding the trivial finite cases for $n=1,2$). The Hantzsche-Wendt manifold corresponds to $F(2,6)=\left\langle a_1,\ldots, a_6\,\middle|\,a_{i}a_{i+1}=a_{i+2} \textrm{ (indices mod $6$) }  \right\rangle$ which can be presented in a balanced manner as depicted above.

\newpage

\section{Nil geometry: Heisenberg manifold}

\begin{figure}[h!]
    \centering
    \includegraphics[width=0.7\textwidth]{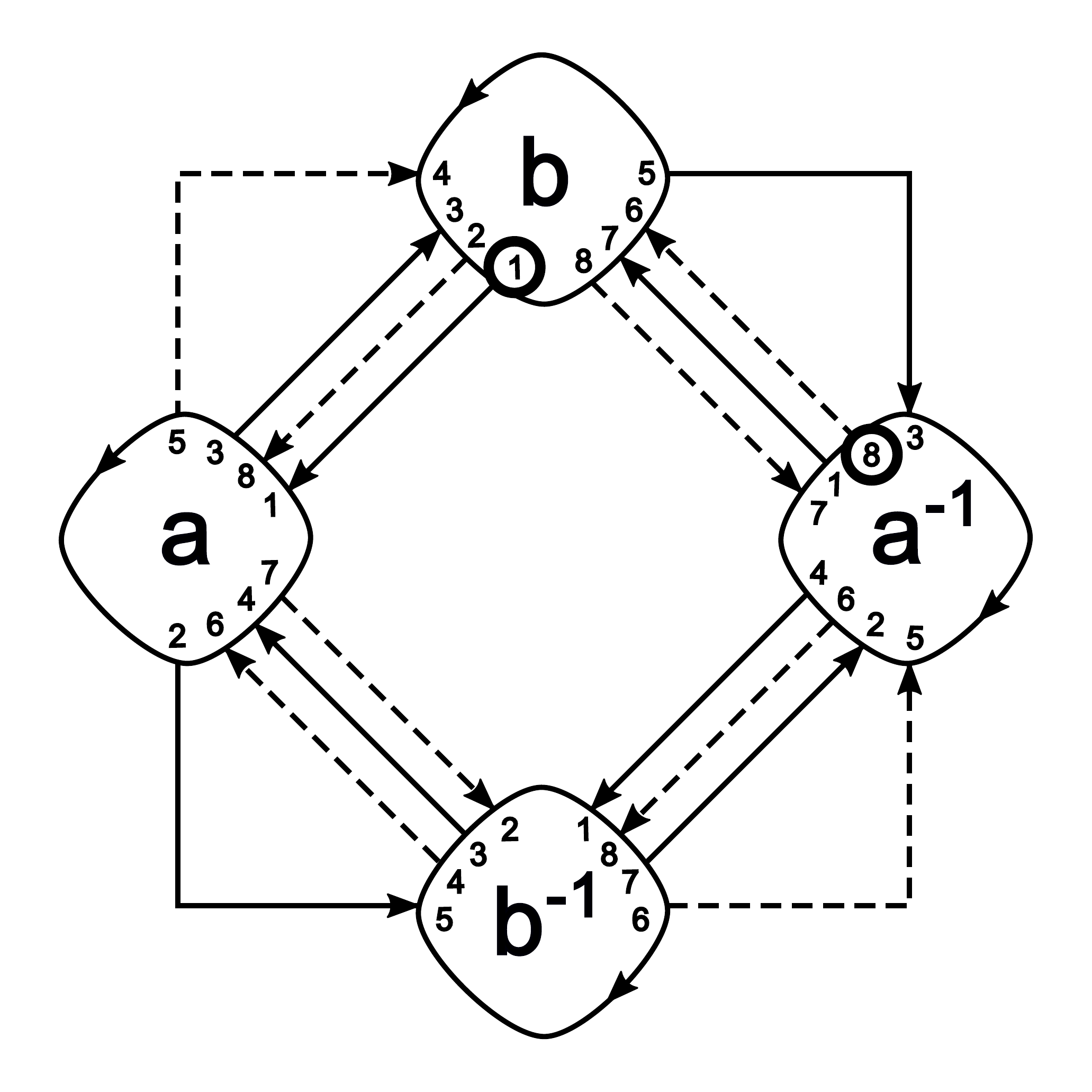}
    \caption{The two relations are $aba^{-1}b^{-1}a^{-1}bab^{-1}=1$ and $ba^{-1}bab^{-1}a^{-1}b^{-1}a=1$.}
    \label{fig:heisenberg}
\end{figure}

The Heisenberg manifold is the quotient of the (unique non-abelian nilpotent) 3-dimensional group of upper triangular $3\times 3$ real matrices by its integer lattice. This lattice has a well known presentation of $\left\langle a,b,c\,\middle|\,c=aba^{-1}b^{-1}, ac=ca, bc=cb\right\rangle$ which simplifies to the one given above.

\newpage

\section{Sol geometry: a mapping torus of a torus}

\begin{figure}[h!]
    \centering
    \includegraphics[width=0.7\textwidth]{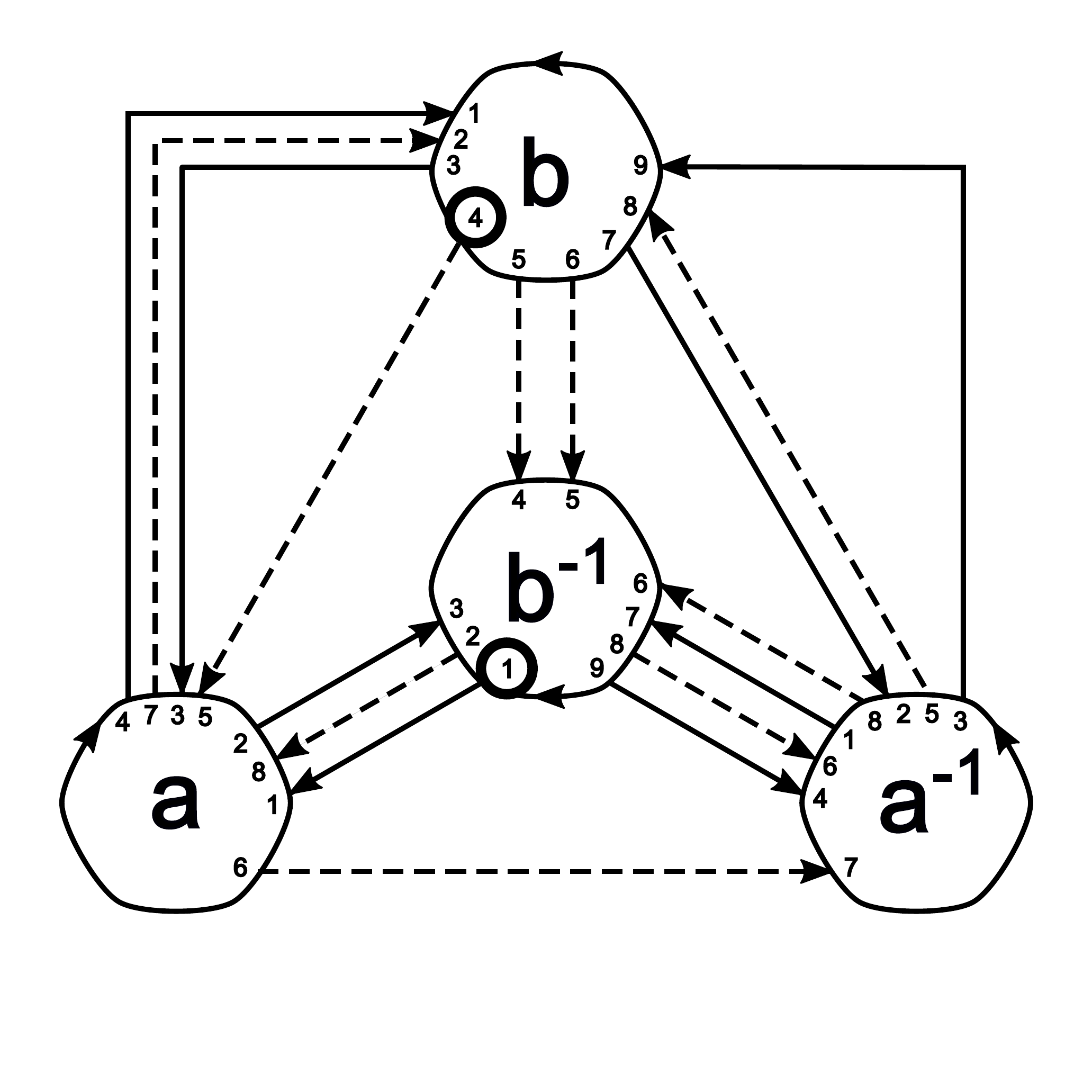}
    \caption{The two relations are $ab^{-1}a^{-1}b^{-1}aba^{-1}b=1$ and $aba^{-2}bab^{-3}=1$.}
    \label{fig:sol}
\end{figure}

Any matrix $L\in Sl(2,\mathbb{R})$ acts by a homeomorphism on a 2-torus, and we can take the mapping torus of this action, $M_L$. By \cite{genusofsol}, Theorem 4.2, $M_L$ having the Sol geometry is equivalent to $L$ being conjugate to $\begin{psmallmatrix} m & -1 \\ 1 & 0 \end{psmallmatrix}$ for some $|m|\geq 3$. Furthermore, for $L=\begin{psmallmatrix}3 & -1\\ 1 & 0\end{psmallmatrix}$ $M_L$ has genus 2 (by Theorem 6.2 in the same paper). Presentations of fundamental groups of torus bundles in terms of the matrix $L$ are given in \cite{presentationofsol} and the one in Section 3 gives in our case $\pi_1(M_L)=\left\langle x,y,t\,\middle|\,xyx^{-1}y^{-1}=1,txt^{-1}=x^3y^{-1},tyt^{-1}=x\right\rangle$, which in turn gives rise to the balanced one we depicted above.

\newpage

\section{\texorpdfstring{$\widetilde{Sl(2,\mathbb{R})}$}{Sl(2,R)} geometry: a Seifert sphere.}

\begin{figure}[h!]
    \centering
    \includegraphics[width=0.7\textwidth]{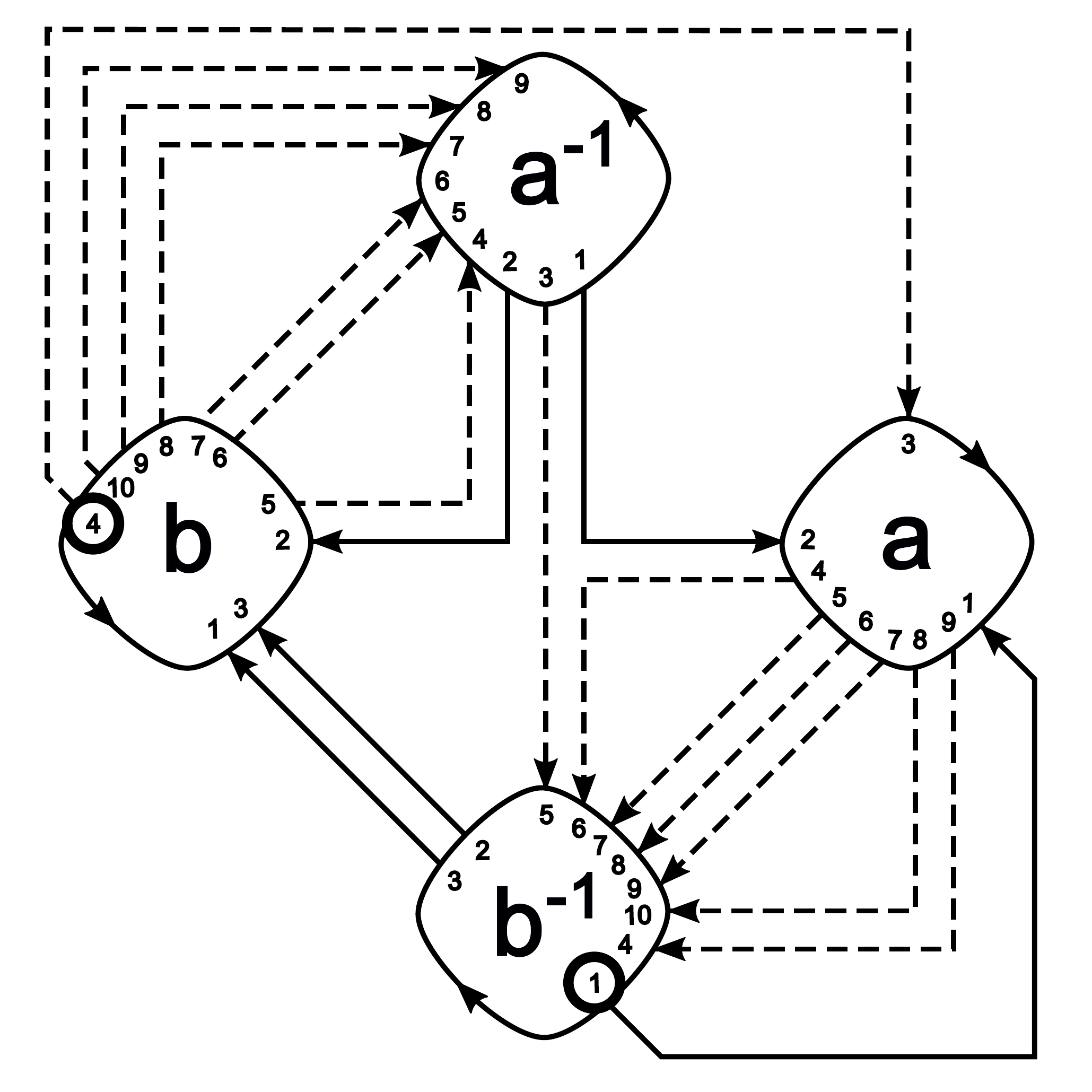}
    \caption{The two relations are $a^2b^3=1$ and $a(b^{-1}a^{-1})^6b^{-1}=1$.}
    \label{fig:brie}
\end{figure}

It is well known that the Brieskorn homology spheres $\Sigma(a_1,a_2,a_3)$, the links in $(0,0,0)$ of $x^{a_1}+y^{a_2}+z^{a_3}=0$ in $\mathbb{C}^3$, admit $\widetilde{Sl(2,\mathbb{R})}$ geometry for coprime $a_i$'s, save for the special cases of $\Sigma(2,3,5)$ yielding the Poincar\'e sphere, and $\Sigma(1,a_2,a_3)$ yielding the standard sphere (and obviously their permuted copies). This is also true in a slightly more general situation of links of $b_1 x^{a_1}+b_2y^{a_2}+b_3z^{a_3}=0$ (sometimes called Seifert spheres), provided that all $b_i$'s are non-zero, and $a_1a_2a_3\left(\frac{b_1}{a_1}+\frac{b_2}{a_2}+\frac{b_3}{a_3}\right)=\pm 1$. One such link, for $x^2-y^3-z^7=0$, turns out to be the Seifert fibration of the symbol $\{0,(o1,0);(2,1),(3,-1),(7,-1)\}$ and as such its fundamental group can be presented in the usual manner (cf. \cite{orlik}) as
$$
\langle a,b,c,d \left|\right. 1=abc=ada^{-1}d^{-1}=bdb^{-1}d^{-1}=cdc^{-1}d^{-1}=a^2d=b^3d^{-1}=c^7d^{-1}\rangle
$$
One can also check directly that this group is perfect, that the orbifold Euler characteristic of this Seifert fibration is negative, and that $\frac{1}{2}-\frac{1}{3}-\frac{1}{7}$ is non-zero, cf. \cite{Scott}, so indeed provides an example of a homology sphere with a $\widetilde{Sl(2,\mathbb{R})}$ geometry. In any case, this presentation is readily simplified to the one depicted.

\newpage

\section{$\mathbb{H}^2\times\mathbb{R}$ geometry: a Seifert fibration over a sphere.}

\begin{figure}[h!]
    \centering
    \includegraphics[width=0.7\textwidth]{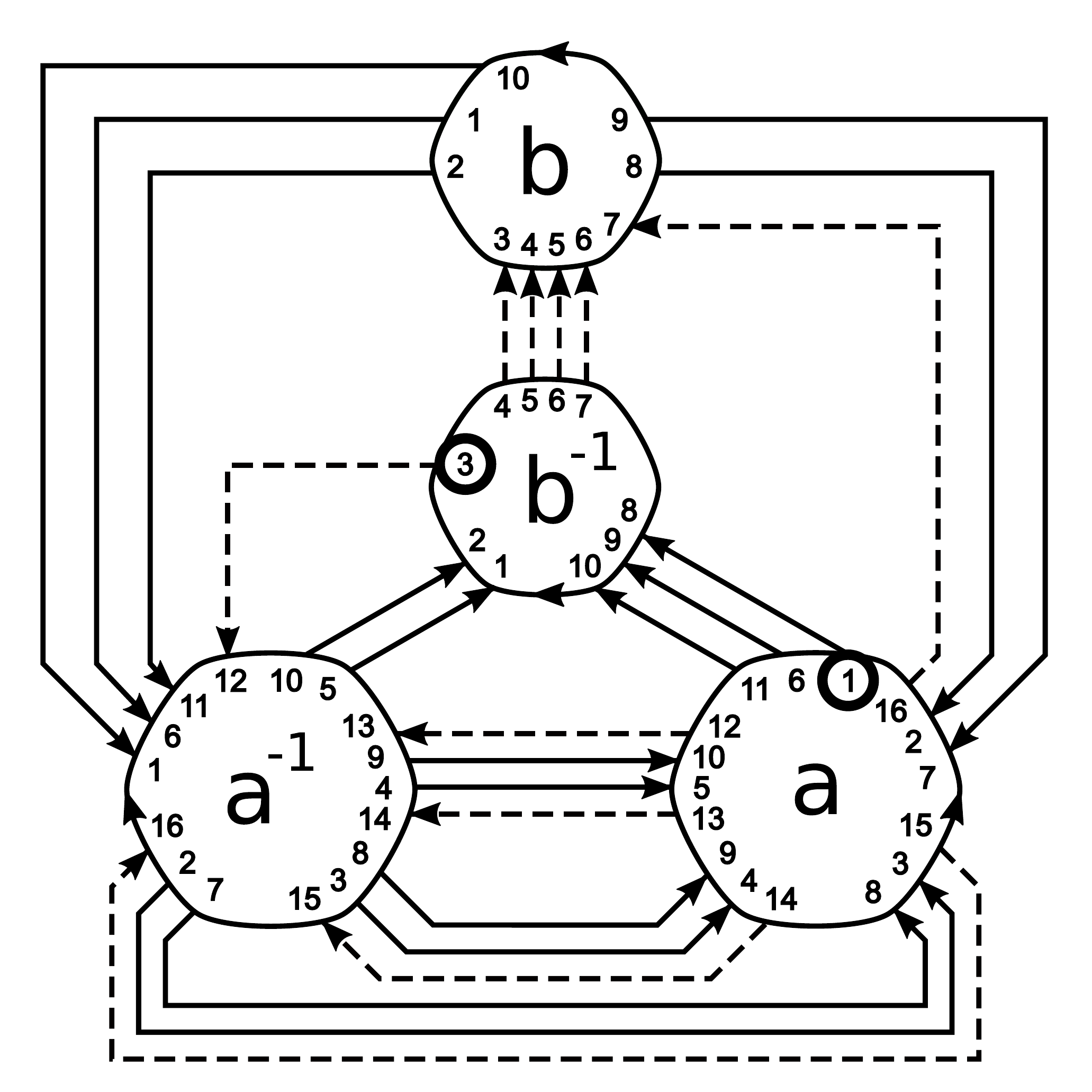}
    \caption{The two relations are $b^{-1}a^4b^{-1}a^{-1}b^{-1}a^4b^{-1}a^{-1}b^{-1}a^{-1}$ and $b^5a^{-5}=1$.}
    \label{fig:seifert}
\end{figure}

Somewhat surprisingly, the most complicated case is that of the $\mathbb{H}^2\times\mathbb{R}$ geometry. By the classical classification of Seifert fibered spaces and their well-known geometric properties (cf. \cite{Scott}), a 3-manifold has this geometry iff it is Seifert fibered, has negative orbifold Euler characteristic, and integer Euler number. By Theorem 1.1 in \cite{Boileau-Zieschang}, every such Seifert fibration over a sphere with exactly three exceptional fibers $(a_i,b_i)$ will have a genus two Heegaard decomposition, and there are countably many such examples (each given by any positive integer solution to $\frac{1}{a_1}+\frac{1}{a_2}+\frac{1}{a_3}<1$ and $\frac{b_1}{a_1}+\frac{b_2}{a_2}+\frac{b_3}{a_3}\in\{1,2\}$, with coprime $b_i<a_i$). In drawing a concrete such diagram, however, the caveat is that the relations in the usual presentation of the fundamental group grow large with $a_i$'s and $b_i$'s, and furthermore these presentations behave stubbornly and seem to resist becoming balanced on two generators under standard algorithms (we base this comment on several experiments with GAP and Heegaard program). We were able to draw one minimal example: depicted above is the Seifert space $\left\{ -1,(o1,0);(5,1),(5,1),(5,3)\right\}$. The standard presentation of its fundamental group is
$$
\langle a,b,c,d \left|\right. 1=abcd=ada^{-1}d^{-1}=bdb^{-1}d^{-1}=cdc^{-1}d^{-1}=a^5d=b^5d=c^5d^3\rangle
$$
which immediately reduces to
$$
\langle a,b \left|\right. 1=ba^5b^{-1}a^{-5}=b^5a^{-5}=b^{-1}a^4b^{-1}a^4b^{-1}a^4b^{-1}a^4b^{-1}a^{-11}\rangle
$$
in which the first relation is a consequence of the second, but we leave it for the Reader's convenience, as it helps to show that the third relation can be rewritten as in the diagram's caption.

\printbibliography

\end{document}